\documentclass[a4paper,12pt]{amsart}
\usepackage{amsmath}
\usepackage{amssymb}
\usepackage{amsfonts}

\DeclareMathSymbol{\subsetneqq}{\mathbin}{AMSb}{36}

\DeclareMathSymbol{\subsetneqq}{\mathbin}{AMSb}{36}

\newtheorem{thm}{{\bf Theorem}}[section]
\newtheorem{lem}{{\bf Lemma}}[section]

\theoremstyle{remark}

\theoremstyle{definition}

\author{Olfa Bjaoui}
\address{University of Tunis ElManar, Faculty of Sciences of Tunis,
Department of Mathematics, ElManar 2092, Tunisia} \email{bjaouiolfa@yahoo.fr}
\thanks{O. B. is grateful to the Laboratory of
PDE and Applications at the Faculty of Sciences of Tunis.}

\author{Mohamed Majdoub}
\address{University of Tunis ElManar, Faculty of Sciences of Tunis,
Department of Mathematics, ElManar 2092, Tunisia} \email{mohamed.majdoub@fst.rnu.tn}
\thanks{M. M. is grateful to the Laboratory of
PDE and Applications at the Faculty of Sciences of Tunis.}

\title[Global weak solutions ...]
{Global weak solutions for some Oldroyd models}
\date{\today}
\begin{document}
\begin{abstract}
 We investigate an evolutive system of non-linear partial differential equations derived from Oldroyd models on Non-Newtonian flows. We prove global existence of weak solutions, in the case of a smooth bounded domain, for general initial data. The results hold true for the periodic case.
\end{abstract}




\maketitle \tableofcontents


\section{Introduction}
The Navier-Stokes equation is a model that describes the evolution  of incompressible fluids with constant viscosity called Newtonian fluids. However, many experiments show that the viscosity of a fluid may vary with the pressure, see Andrade \cite{Andrade}, see also the book by Bridgman \cite{Bridgman}. Further details and references to more recent experimental studies can be found in the book by Szeri \cite{Szeri} and in the paper by Malek and Rajagopal \cite{Mal+Raja}. There are also many experiments which show that the viscosity may depend on the symmetric part of the velocity gradient. We can refer to Schowalter \cite{Scho}, Huilgol \cite{Huil}. Recently, Malek, Necas and Rajagopal \cite{Mal+Raja+Necas}, Bulicek, Majdoub and Malek \cite{Mal+Maj+Buli} have established existence results concerning the flows of fluids called non-Newtonian fluids, whose viscosity depends on both the pressure and the symmetric part of the velocity gradient. The constitutive equation for such fluids is given by
\begin{equation}\label{eq con}T=-p I + \nu (p, |D(v)|^{2})D(v),\end{equation}
where $T$ is the Cauchy stress, $p$ is the pressure, $v$ is the velocity, $D(v)$ is the symmetric part of the velocity gradient
and $\nu$ is the viscosity function. Thus, the evolution of such fluids is governed by the equation$$\rho(\partial_t v+v.\nabla v)=div T,$$where $\rho$ is the density of the fluid.\\However, there are some fluids that do not obey to the constitutive equation \eqref{eq con}
such as Blood. Yet, it was shown experimentally that Blood is a complex rheological mixture that exhibits shear thinning and elastic behavior, see Thurston \cite{Thurston}. The constitutive equation for such fluids is given by $$T=-p I +S,$$where $S$ is the extra-stress tensor which is related to the kinematic variables through
$$S+\lambda_1 \frac{DS}{Dt}=2\mu (|D(v)|^2)D(v)+2 \lambda_2\frac{DDv}{Dt},$$
$\mu$ is the viscosity function, $\lambda_1>0$ and $\lambda_2>0$ are viscoelastic constants.\\The symbol $\frac{D}{Dt}$ denotes the objective derivative of Oldroyd type defined by$$\frac{D S}{Dt}=\partial_t S+v.\nabla S+S.w(v)-w(v)S,$$where $w(v)$ is the antisymmetric part of the velocity gradient. The extra-stress tensor $S$ is decomposed into the sum of its Newtonian part $$\tau_s=2 \frac{\lambda_2}{\lambda_1}D(v)$$ and its viscoelastic part $\tau_e$. The constitutive equation for $\tau_e$ is given by $$\tau_e+\lambda_1 \frac{D \tau_e}{Dt}=2 \Big(\mu (|D(v)|^2)-\frac{\lambda_2}{\lambda_1}\Big) D(v).$$
Hence, we get the generalized Oldroyd-B model given by
\begin{eqnarray*}
\begin{cases}
\rho\partial_t v+ \rho v.\nabla v- \frac{\lambda_2}{\lambda_1}\Delta
v+\nabla
p=div(\tau_e),\\
\partial_t  \tau_e+ \frac{1}{\lambda_1}\tau_e+v.\nabla \tau_e+\tau_e.w(v)-w(v)\tau_e=\frac{2}{\lambda_1}\Big(\mu (|D(v)|^2)-\frac{\lambda_2}{\lambda_1}\Big)D(v).
\end{cases}
\end{eqnarray*}
Particularly, when $\mu$ is constant, we recover the Oldroyd-B fluid with constant viscosity, see \cite{Lions+Masm}. Existence of local strong solutions to the Oldroyd-B model was proved by Guillop\'e and Saut in \cite{Gui+Saut} and \cite{Gui+Saut'}. Fernandez-Cara, Guillen and Ortega \cite{Fercara+Gui+ort}, \cite{Fercara+Gui+ort'} and \cite{Fercara+Gui+ort''} proved local well posedness in Sobolev spaces. In the frame of critical Besov spaces, Chemin and Masmoudi had proved local and global well-posedness results in \cite{chem+masm}. In addition, non-blow up criteria for Oldroyd-B model were given in \cite{kupfer+m+t} as well as in \cite{Lei+mas+zhou}. For the sake of completeness, we refer the reader to the most recent papers about Oldroyd models as \cite{Lei, Lei-Arxiv, Lei+Liu+Zhou', Lei', Lei+Wang}.\\Global existence for small data was proved in \cite{Lei+Liu+Zhou} and \cite{Lei+Zhou}. Considering general initial data, Lions and Masmoudi \cite{Lions+Masm} had established results of global existence of weak solutions when the viscosity function $\mu$ is constant.
Our aim in this paper is to generalize the results in \cite{Lions+Masm}. More precisely, we will focus on the following system
\begin{eqnarray*}(S)
\begin{cases}
\partial_t v+v.\nabla v-\textmd{div }(f(D(v)))
+\nabla p=\textmd{div }\tau,\,\,\,
\mbox{in}\,\,\,\mathbf{R}^{+}\times\Omega,\\
\partial_t \tau+v.\nabla \tau+\tau.w(v)-w(v).\tau+a \tau=g(D(v)),\,\,\,
\mbox{in}\,\,\,\mathbf{R}^{+}\times\Omega,\\
\textmd{div }{v}=0,\,\,\,\mbox{in}\,\,\, \mathbf{R}^+\times\Omega,
\\v(0,x)=v_0(x),\, \tau(0,x)=\tau_0(x),\end{cases}
\end{eqnarray*}
where $\Omega$ can be considered either the torus $\mathrm{T}^{n}$, or a smooth bounded domain of $\mathrm{R}^{n}$, $n=2,\, 3$ and in this case $(S)$ is supplemented by the Dirichlet homogeneous boundary condition.\\
The function $g(D(v))$ is given by
$$g(D(v))=\tilde{\mu}(|D(v)|^{2})D(v)=\frac{b}{1-\theta}(\mu(|D(v)|^{2})-\theta)D(v),$$
where $\mu(|D(v)|^{2})=1-\lambda+\lambda (1+|D(v)|^{2})^{\frac{r-2}{2}},$ $\lambda \in [0,1],r\in [1,2]$.\\
Physically, the parameters $a$ and $b$ are given by $a=\frac{1}{We}, b=\frac{2(1-\theta)}{We}$, where $We$ is the Weissenberg number and $\theta \in ]0,1[$ is the ratio between the so called relaxation and times.\\
Let us mention that the expression of $\tilde{\mu}$ is motivated by the system studied by Arada and Sequeira \cite{ar+sq} in the steady case with $f$ being the identity map, where existence of a unique solution was established for small and suitably regular data. Notice that when the term $f(D(v))$ is replaced by $\nu \Delta v$, with $\nu= \frac{\theta}{Re}$, $Re$ being the the Reynolds number, and $r$ is equal to $2$, then $(S)$ turns into the system studied by Lions and Masmoudi in \cite{Lions+Masm}. Let us remark that when $b=0$ and $\tau_0=0$, then $\tau=0$ solves the second equation in $(S)$
which consequently will be reduced to the system studied in \cite{Drey+Hung} and \cite{Drey+Hung'}.\\Our objective is to prove global existence of weak solutions, for general initial data, under
suitable hypotheses on the function $f$.\\The layout of this paper is as follows. In the following section, we give some notations and we introduce the functional
spaces used along this paper. The third section is devoted to some technical lemmas and to the statement of the main results. In the fourth section, we prove
existence of approximate solutions $(v^N,\tau^N)$ to the system $(S)$ by using Galerkin method. Then, we prove the equi-integrability in the Lebesgue space $L^{2}([0,T]\times \Omega)$ of the
sequence $\tau^N$, and we derive an evolution equation for $\eta:=\overline{|\tau^N -\tau|^2}$. This enables us to get the strong convergence of the sequence $\tau^N$ in $L^{2}([0,T]\times \Omega)$. Therefore, we can identify the weak limits in the sense of distributions
of non linear terms.
\section{Notations and Functional spaces}
Let us introduce some notations which we use throughout this paper. With $M^{n\times n}$ we denote the real vector space of $n\times n$ matrices.
We have $(div \tau)_{i}=\sum_{j}\partial_{j}\tau_{ij}$, for $\tau \in M^{n\times n}$, and for a vector field $v$, we have $(\nabla u)_{ij}=\partial_{j}u_{i}$.\\
For $A=(a_{i,j})$ and $B=(b_{i,j})$ in $M^{n\times n}$, $A:B$ stands for the sum $\sum_{i,j}a_{i,j}b_{i,j}$ and $|A|^2$ denotes the quantity $A:A$.\\Let us denote by:$$\mathcal{V}=\{\phi \in \mathcal{D}(\Omega)^{n}, div \phi =0\},$$
and$$\widetilde{\mathcal{V}}=\{M \in \mathcal{D}(\Omega)^{n\times n},  M^t=M \},$$ when $\Omega$ is a bounded domain of $\mathbf{R}^{n}$, $$\mathcal{V}=\{\phi \in \mathbf{C}^{\infty}_{per}(\Omega)^{n}, div \phi =0, \int_{\Omega} \phi dx=0\},$$
and$$\widetilde{\mathcal{V}}=\{M \in \mathbf{C}^{\infty}_{per}(\Omega) ^{n \times n},  M^t=M, \int_{\Omega} M dx=0\}$$
when $\Omega$ is the torus $T^{n}$,\\$\mathcal{H}=$ the closure of $\mathcal{V}$ in the $L^{2}(\Omega)^{n}$-norm,\\
$\widetilde{\mathcal{H}}=$ the closure of $\widetilde{\mathcal{V}}$ in the $L^{2}(\Omega)^{n \times n}$-norm,\\$ \textbf{V}_q=$ closure of $\mathcal{V}$ in the $L^q(\Omega)^{n^2}$-norm of gradients, $q\geq 1$.\\
$\widetilde {\textbf{V}}_{q}=$ closure of $\widetilde{\mathcal{V}}$ in the $L^q(\Omega)^{n\times n^2}$-norm of gradients, $q\geq 1$.\\
$V_s=$ the closure of $\mathcal{V}$ with respect to the $W^{^{s,2}}(\Omega)^{n}$-norm, with $s>1+\frac{n}{2}$,\\
$\widetilde{V}_s=$ the closure of $\widetilde{\mathcal{V}}$ with respect to the $W^{^{s,2}}(\Omega)^{n\times n}$-norm, with $s>1+\frac{n}{2}$.\\
The condition on $s$ is due to the fact that: if $v \in W^{s,2}(\Omega)^{n}$, then $\nabla v \in W^{s-1,2}(\Omega)^{n^2}$ and $W^{s-1,2}(\Omega) \hookrightarrow L^\infty(\Omega)$ if $\frac{1}{2}-\frac{s-1}{n}< 0$.\\
We denote by $V_s'$ the dual space of $V_s$ and by $<,>_{s}$ the duality between $V_{s}$ and $V_s'$.
The scalar product in $L^{2}(\Omega)$ will be denoted by $(.,.)$.\\
For $s>1+\frac{n}{2}$ and $q\geq 2$, we have the following inclusions$$V_s \subset \textbf{V}_q \subset \mathcal{H} \simeq \mathcal{H}' \subset \textbf{V}'_q \subset V'_s,$$
and $$\tilde{V}_s \subset \widetilde {\textbf{V}}_{q} \subset \widetilde{\mathcal{H}} \simeq \widetilde{\mathcal{H}}' \subset \widetilde {\textbf{V}}_{q}' \subset \tilde{V}_s'.$$
For a sequence $f_N$ in $\mathcal{D}'([0,T]\times \Omega)$, we denote by $\overline{{f_N}}$ its weak limit in the sense of distributions.\\
By $C$ we denote any constant that may depend on $|\Omega|$ the n-dimensional Lebesgue measure of $\Omega$, $v_0$, $ \tau_0$ and on $T> 0$, but not on $N$.\\
Let us mention that subsequences will not be relabeled.
\section{Technical lemmas and statement of the main results}
We give some technical lemmas needed for the proofs of our main results.\\
One of the properties of the norm on Banach spaces is the lower semi-continuity given by the following lemma. For more details, we refer to \cite{Yosida}.
\begin{lem}\label{lim inf}
Let $X$ be a Banach space equipped with a norm $\|.\|$. Let $(x_N)$ be a sequence in $X$ that converges weakly to some $x$ in $X$, then $(x_N)$ is bounded in $X$ and we have $$\|x\|\leq \liminf _{N \rightarrow +\infty} \|x_N\|.$$\end{lem}
\begin{lem}\label{Vitali}(Vitali's Lemma)
Let $\Omega$ be a bounded domain in $\mathbf{R}^{n}$ and $f_N:\Omega \rightarrow \mathbf{R}$ be integrable for every $N \in \mathbf{N}$. Assume that\\
(i)\quad $\lim_{N \rightarrow +\infty}f_N(y)$ exists and is finite for almost all $y \in \Omega$,\\
(ii)\quad for every $\varepsilon > 0$, there exists $\delta > 0$ such that$$\sup _{N \in \mathbf{N}}\int_H |f_N(y)| dy <\varepsilon \quad \quad \forall H \subset \Omega, |H| < \delta.$$Then$$\lim_{N \rightarrow +\infty}\int_{\Omega}f_N(y)dy= \int_{\Omega}\lim_{N \rightarrow +\infty}f_N(y)dy.$$
\end{lem}
For the proof of lemma \ref{Vitali}, we can refer to \cite{Alt}.
\begin{lem}\label{Aubin-Lions}(Aubin-Lions lemma) Let $1<\alpha,\beta  <+\infty$ and $T> 0$. Let $X$ be a Banach space, and let $X_0, X_1$ be separable and reflexive Banach spaces such that $X_0$ is compactly embedded into $X$ which is continuously embedded into $X_1$, then$$\{v\in L^{\alpha}([0,T], X_0); \partial_t v \in  L^{\beta}([0,T], X_1)\} \,\,\mbox{is compactly embedded into} \,\,L^{\alpha}([0,T], X).$$
\end{lem}
For the proof of lemma \ref{Aubin-Lions}, we refer to \cite{Lions 69}. A generalized form of this lemma for locally convex spaces and $\beta=1$ can be found in \cite{Rubicek}.\\
The following lemma plays an important role in the theory of existence of solutions to ordinary differential equations. The proof of such lemma can be found in
\cite{Walter}.
\begin{lem}\label{Caratheodory }(Caratheodory Theorem)
Let $c:I_{\delta}\equiv [t_0-\delta,t_0+\delta]\rightarrow \mathbf{R}^{n}$ the system of ordinary differential equations
\begin{eqnarray*}(E)
\begin{cases}

\frac{d}{dt}c(t)=F(t,c(t)), \,\,\,t\in I_{\delta}\,\,\,\\
c(t_0)=c_0 \in \mathbf{R}^{n},
\end{cases}
\end{eqnarray*}
Assume $F:I_{\delta} \times K \rightarrow \mathbf{R}^{n}$, where $K\equiv \{c\in \mathbf{R}^{n}, |c-c_0|< \beta\}$, for some $\beta>0$.\\
If $F$ satisfies the Caratheodory conditions:\\
(i)$\quad t\mapsto F_i(t,c)$ is measurable for all $i=1,...,n$ and for all $c \in K$,\\
(ii)$\quad c\mapsto F(t,c)$ is continuous for almost all $t\in I_{\delta}$,\\
(iii)\quad there exists an integrable function $G:I_{\delta}\rightarrow \mathbf{R}$ such that $$|F_i(t,c)| \leq G(t), \quad \forall (t,c)\in I_{\delta} \times K, \quad \forall i=1,...,n,$$
then there exists $\delta' \in ]0, \delta[$ and a continuous function $c:I_{\delta'}\rightarrow \mathbf{R}^{n}$ such that\\
(i)$\quad \frac{dc}{dt}$ exists for almost all $t\in I_{\delta'}$,\\(ii)$\quad c$ solves $E$.
\end{lem}
The following theorem, the proof of which can be found in  \cite{Bresis}, deals with compact injections of Sobolev spaces into Lebesgue ones.
\begin{thm}\label{Rellich}(Rellich-Kondrachov) Let $\Omega$ be a $C^{1}$ bounded domain of $\mathbf{R}^{n}$. We have the following compact injections\\
(i)\,\,if $p=n$, then $W^{1,p}(\Omega)\subset L^q(\Omega)$, $\forall q \in [1,+\infty[$,\\(ii)\,if $p >n$, then $W^{1,p}(\Omega)\subset C(\overline{\Omega})$.
\end{thm}
In order to prove the existence of approximate solutions, we use Galerkin method by considering a special basis in the space $V_{s}(\Omega)$ as in \cite{Drey+Hung}. The proof of existence of such basis, see \cite{ma+ne+ro+ru},
relies on solving the following spectral problem: find $w^r \in V_{s}$ and $\lambda_r \in \mathbf{R}$ satisfying $$<w^r,\psi>_{s}=\lambda_r(w^r,\psi), \quad \forall \psi \in V_s.\quad (\mathcal{P})$$
\begin{thm}\label{basis}
There exists a countable set $\{\lambda_r\}_{r=1}^{\infty}$ and a corresponding family of eigenvectors $\{w^r\}_{r=1}^{\infty}$ solving the problem $(\mathcal{P})$ such that\\
(i)$\quad (w^r,w^{r'})=\delta_{r,r'}\quad \forall r,r' \in \mathbf{N}$,\\
(ii)$\quad 1\leq \lambda_1 \leq \lambda_2\leq ...$ and $\lambda_r \rightarrow \infty$ as $r\rightarrow \infty$,\\
(iii)$\quad <\frac{w^r}{\sqrt{\lambda_r}},\frac{w^{r'}}{\sqrt{\lambda_{r'}}}>_{s}=\delta_{r,r'} \quad \forall r,r' \in \mathbf{N}$,\\
(iv)$\quad \{w^r\}_{r=1}^{\infty}$ forms a basis in $V_s$.\\
Moreover defining $\mathcal{H}^N :=span \{w^1,...,w^N\}$ (a linear hull) and $$P^N(v):=\sum_{i=1}^{N}(v,w^r)w^r \in \mathcal{H}^N,\, v \in V_{s},$$ we obtain
$$\|P^N\|_{\mathcal{L}(V_s, V_s)}\leq1,\, \|P^N\|_{\mathcal{L}(V'_s, V'_s)}\leq1,\, \|P^N\|_{\mathcal{L}(\mathcal{H}, \mathcal{H})}\leq1.$$
\end{thm}
The same results hold if we replace $V_s$ by $\widetilde{V}_s$.\\
One of the tools that we use in order to prove our main results is Young measures. Thus, we should give some well-known facts about this tool. We can refer to
\cite{young2} as well as to \cite{young1} and \cite{ma+ne+ro+ru} for more details.
\begin{thm}\label{Dun Pet}
Let $(f_N:\Omega \rightarrow \mathbf{R}^n)$ be a bounded sequence in $L^1(\Omega)$ and $\overline{f_N}$ be its weak limit in $\mathcal{D}'(\Omega)$. Assume that  $(f_N)$  is equi-integrable in $L^1(\Omega)$, then $$\overline{f_N}=\int_{\mathbf{R}^n}\lambda d\nu_x(\lambda),\quad a.e \,x \in \Omega,$$where $\nu$ is the Young measure generated by the sequence $f_N$.\\Now, assume that $f_N$ converges weakly in $L^1(\Omega)$, then $f_N$ is equi-integrable in $L^1(\Omega)$ and
$$\overline{f_N}=\int_{\mathbf{R}^n}\lambda d\nu_x(\lambda),\quad a.e \,x \in \Omega.$$
\end{thm}
\begin{thm}\label{fatou}
Let $(f_N:\Omega \rightarrow \mathbf{R}^n)$ be a sequence of maps that generates the Young measure $\nu$. Let $F: \Omega \times \mathbf{R}^n\rightarrow \mathbf{R} $ be a Caratheodory function that is a function that satisfies the Caratheodory conditions. Assume that the negative part $F^-(.,f_N(.))$ is weakly relatively compact in $L^1(\Omega)$, then$$\int_{\Omega}\int_{\mathbf{R}^n}F(x,\lambda)d\nu_x(\lambda)dx \leq \liminf_{N \rightarrow +\infty}\int_{\Omega} F(x,f_N(x))dx.$$
\end{thm}
Now, we are ready to state our main results. The first Theorem gives a global existence result for the following system
\begin{eqnarray*}(S_1)
\begin{cases}
\partial_t v+v.\nabla v-\textmd{div }(f(D(v)))
+\nabla p=\textmd{div }\tau,\,\,\,
\mbox{in}\,\,\,\mathbf{R}^{+}\times\Omega,\\
\partial_t \tau+v.\nabla \tau+a \tau=g(D(v)),\,\,\,
\mbox{in}\,\,\,\mathbf{R}^{+}\times\Omega,\\
\textmd{div }{v}=0,\,\,\,\mbox{in}\,\,\, \mathbf{R}^+\times\Omega,
\\v(0,x)=v_0(x),\, \tau(0,x)=\tau_0(x).\end{cases}
\end{eqnarray*}
\begin{thm}\label{global1}
Let $f:M^{n\times n}\rightarrow M^{n\times n}$ be a continuous function satisfying the following hypotheses for some $p\in ]2,+\infty[$ when $n=2$ and $p\in ]\frac{5}{2},+\infty[$ when $n=3$:\\
$(H_{1})$ growth: there exists $c> 0$ and $\tilde{c}\geq 0$ such that $$|f(A)|\leq \tilde{c}+c |A|^{p-1},\,\forall A \in M^{n\times n},\, f(0)=0,$$
$(H_{2})$ monotonicity: there exists $\nu> 0$ such that$$(f(A)-f(B)):(A-B)\geq \nu (|A-B|^{2}+|A-B|^{p}),\,\forall A,\, B\in M^{n\times n}.$$
Let $v_{0}$ and $\tau_{0}$ be in $\mathcal{H}(\Omega)$ and $\widetilde{\mathcal{H}}(\Omega)$ respectively.\\
i)If $\nu$ satisfies $2 \nu (1-\theta)> 1$, then for an arbitrary $\lambda \in [0,1]$, there exists a global weak solution $(v,\tau)$ to the system $(S_1)$ such that $$v \in L^{\infty}(R^{+},\mathcal{H}(\Omega))\cap  L^{p}(0,T,\textbf{V}_p(\Omega)), \,\,\tau \in L^{\infty}(R^{+},\widetilde{\mathcal{H}}(\Omega)), \, \forall T > 0.$$
ii)If $\nu$ is such that $0\leq 2 \nu (1-\theta) \leq 1$, then the same result holds for an arbitrary $\lambda \in [0,\sqrt{2 \nu (1-\theta)}[$.
\end{thm}
Let us assume the existence of a scalar function $U \in C^2 (\mathbf{R}^{n^2})$, called potential of $f$, such that for some $p> 1$, $C_1>0$, $C_2>0$ we have for all $\eta, \xi \in M^{n\times n}_{sym}$ and $i,j,k,l \in \{1,...,n\}$
\begin{equation}\label{eq t1}
\frac{\partial U(\eta)}{\partial \eta_{ij}}=f_{ij}(\eta)
,\end{equation}
\begin{equation}\label{eq t2}
U(0)=\frac{\partial U(0)}{\partial \eta_{ij}}=0
,\end{equation}
\begin{equation}\label{eq t3}
\frac{\partial^2 U(\eta)}{\partial\eta_{ij}\partial\eta_{kl}}\xi_{ij}\xi_{kl}\geq C_1 (1+|\eta|)^{p-2}|\xi|^2,
\end{equation}
\begin{equation}\label{eq t4}
\Big|\frac{\partial^2 U(\eta)}{\partial\eta_{ij}\partial\eta_{kl}}\Big|\leq C_{2}(1+|\eta|)^{p-2}.
\end{equation}
As consequence of these assumptions, there exists $C > 0$ such that
$$|f(A)|\leq C(1+ |A|)^{p-1},$$
Moreover, for $p \geq 2$, there exists $\nu > 0$ such that
$$(f(A)-f(B)):(A-B)\geq \nu (|A-B|^{2}+|A-B|^{p}),\,\forall A,\, B\in M^{n\times n}.$$
Standard examples of functions $f$ whose potentials satisfy these assumptions are
$$f(A)=(1+ |A|)^{p-2}A$$ and $$f(A)=(1+ |A|^2)^{\frac{p-2}{2}}A.$$
For more details about the existence of a such potential $U$ and consequences of properties \ref{eq t1}-\ref{eq t4} we refer to \cite{ma+ne+ro+ru}.\\
The second Theorem concerns the following system
\begin{eqnarray*}(S_2)
\begin{cases}
\partial_t v+v.\nabla v-\textmd{div }(f(D(v)))
+\nabla p=\textmd{div }\tau,\,\,\,
\mbox{in}\,\,\,\mathbf{R}^{+}\times\Omega,\\
\partial_t \tau+v.\nabla \tau+a \tau=bD(v),\,\,\,
\mbox{in}\,\,\,\mathbf{R}^{+}\times\Omega,\\
\textmd{div }{v}=0,\,\,\,\mbox{in}\,\,\, \mathbf{R}^+\times\Omega,
\\v(0,x)=v_0(x),\, \tau(0,x)=\tau_0(x).\end{cases}
\end{eqnarray*}
\begin{thm}\label{global2}
Let $f:M^{n\times n}\rightarrow M^{n\times n}$ be a $C^1$-function satisfying the following hypotheses for some $p\in ]2,+\infty[$ when $n=2$ and $p\in ]\frac{5}{2},+\infty[$ when $n=3$:\\
$(H_{1})$ growth: there exists $c> 0$ such that $$|f(A)|\leq c |A|^{p-1},\quad\forall A \in M^{n\times n},$$
$(H_{2})$ coercivity: there exists $\nu > 0$ such that $$f(A):A\geq \nu |A|^{p},\,\forall A \in M^{n\times n},$$
$(H_{3})$ monotonicity $$(f(A)-f(B)):(A-B)\geq 0,\,\forall A,\, B\in M^{n\times n}.$$
Let $v_{0}$ and $\tau_{0}$ be in $\mathcal{H}(\Omega)$ and $\widetilde{\mathcal{H}}(\Omega)$ respectively. Then, there exists a global weak solution $(v,\tau)$ to the system $(S_2)$ such that $$v \in L^{\infty}(R^{+},\mathcal{H}(\Omega))\cap  L^{p}(0,T,\textbf{V}_p(\Omega)), \,\,\tau \in L^{\infty}(R^{+},\widetilde{\mathcal{H}}(\Omega)), \, \forall T > 0.$$
\end{thm}
Obviously, if $f$ satisfies the hypotheses of Theorem \ref{global1} with $\tilde{c}=0$, then $f$ satisfies the hypotheses of Theorem \ref{global2}.\\Let us remark that the quadratic term $\tau.w(v^N)-w(v^N).\tau$ is not present in $(S_1)$ neither in $(S_2)$. The difficulty of this fact will be explained in the proofs
of Theorems \ref{global1} and \ref{global2} that will be given in the case of a bounded domain. However, they can be easily adapted to the periodic case.
\section{Proofs of the main results}
The results in subsections \ref{Galerkin}, \ref{Uniform}, \ref{Strong1} and \ref{Equi} hold for the system $(S)$ and thus for the systems $(S_1)$ and $(S_2)$.
\subsection{Galerkin approximation}\label{Galerkin}
We will show the existence of approximate solution to the system $(S)$ via Galerkin approximation as in \cite{Drey+Hung}. Hence, let
$\{a^r\}_{r=1}^{\infty}$ and $\{\alpha^r\}_{r=1}^{\infty}$ be basis of $V_{s}$ and $\tilde{V}_{s}$ respectively given by Theorem \ref{basis}.\\
Let $T> 0$ and $N\geq 1$ be fixed. We define $$v^{N}(t,x)=\sum_{k=1}^{N}d_k^N (t)a^k(x),\quad \tau^{N}(t,x)=\sum_{k=1}^{N}c_k^N(t) \alpha^k(x),$$ where the coefficients $c_k^N(t)$ and $d_k^N (t)$ solve the so-called Galerkin system: a system of $2N$ nonlinear equations with $2N$ unknowns
\begin{eqnarray*}(S_N)
\begin{cases}
\int_{\Omega}\partial_t v^N a^i dx+\int_{\Omega}(v^N.\nabla) v^N a^{i}dx+\int_{\Omega}f(D(v^N)):Da^{i} dx
=-\int_{\Omega}\tau^N:D a^i dx,\\
\int_{\Omega}\partial_t \tau^N: \alpha^{j}dx+\int_{\Omega}v^N.\nabla \tau^N :\alpha^{j}dx +a \int_{\Omega}\tau^N : \alpha^{j}dx =\int_{\Omega}\tilde{\mu}(|D(v^N)|^{2})D(v^N):\alpha^{j}dx,
\\c_j^N(0)=(\tau_{0},\alpha^j),d_i^N(0)=(v_{0},a^i), 1\leq i,j\leq N\, \end{cases}
\end{eqnarray*}The initial conditions on $c_k^N$ and $d_k^N$ are such that $$v^N(0,x)=P^N v_0(x), \,\tau^N(0,x)=\widetilde{P}^N\tau_0(x),$$ where $P^N$ and $\tilde{P}^N$ are the orthogonal continuous projectors of $\mathcal{H}$ and $\mathcal{\widetilde{H}}$ respectively onto the linear hulls of the first eigenvectors $a^r, r=1,...N$ and $\alpha^r, r=1,...N$ respectively.\\The orthogonality of the two basis $\{a^r\}_{r=1}^{\infty}$ and $\{\alpha^r\}_{r=1}^{\infty}$ in $\mathcal{H}$ and $\widetilde{H}$ respectively imply that the system $(S_N)$ can be rewritten as
\begin{eqnarray*}
\begin{cases}
\frac{d}{dt}(c_{j}^N,d_i^N) =\mathcal{F}_{i,j}(t,c_1^N,...,c_N^N,d_1^N,...,d_N^N)
\\c_{j}^N(0)=(\tau_{0},\alpha^{j}),d_i^N(0)=(v_{0},a^i),\end{cases}
\end{eqnarray*}
where for $1\leq i,j\leq N$\begin{eqnarray*}\mathcal{F}_{i,j}(t,c_1^N,...,c_N^N,d_1^N,...,d_N^N)&=&\Big(-\int_{\Omega}(\sum_{k=1}^{N}c_k^N \alpha^k):Da^i dx-\int_{\Omega}(\sum_{k=1}^{N}d_k^N a^k). (\sum_{k=1}^{N}d_k^N \nabla a^k)a^i dx\\ &-& \int_{\Omega}f(\sum_{k=1}^{N}d_k^N D(a^k)):D(a^i)dx,
-\int_{\Omega}  (\sum_{k=1}^{N}d_k^N a^k). (\sum_{k=1}^{N}c_k^N \nabla \alpha^k):\alpha^{j}dx\\&-& a \int_{\Omega}(\sum_{k=1}^{N}c_k^N \alpha^k):  \alpha^{j}dx +\int_{\Omega}\tilde{\mu}(|D(\sum_{k=1}^{N}d_k^N a^k)|^{2})D(\sum_{k=1}^{N}d_k^N a^k):\alpha^{j}dx\Big)\\ &+& \int_{\Omega}(\sum_{k=1}^{N}d_k^N w(a^k)).(\sum_{k=1}^{N}c_k^N \alpha^k):\alpha^{j}dx\\&-& \int_{\Omega}(\sum_{k=1}^{N}c_k^N \alpha^k).(\sum_{k=1}^{N}d_k^N w(a^k)):\alpha^{j}dx.\end{eqnarray*}
Let $R> 0$ and $K\subset \mathbf{R}^{2 N}$ be the ball of center $(c_1^N(0),...,c_N^N(0),d_1^N(0),...,d_N^N(0))$ and of radius $R$. We consider $\mathcal{F}_{i,j}:[0,T]\times K \rightarrow \mathbf{R}^{2 N}$.\\The continuity of $\tilde{\mu}$ and $f$ lead to the continuity of $\mathcal{F}_{i,j}$ over $[0,T]\times K$. In addition, thanks to the continuous inclusion $W^{s-1,2}(\Omega) \hookrightarrow L^\infty(\Omega)$ if $\frac{1}{2}-\frac{s-1}{n}< 0$, we get
\begin{eqnarray*}\Big|\int_{\Omega}(\sum_{k=1}^{N}c_k^N \alpha^k):Da^i dx\Big| \leq C(R,N)\|Da^i\|_{L^{\infty}(\Omega)}\sum_{k=1}^{N}\|\alpha^k\|_{L^{2}(\Omega)}. \end{eqnarray*}
In the same way, we have
\begin{eqnarray*}
\Big|\int_{\Omega}(\sum_{k=1}^{N}d_k^N a^k). (\sum_{k=1}^{N}d_k^N \nabla a^k)a^i dx\Big|\leq C(R,N) \|a^i\|_{L^{2}(\Omega)}\sum_{k,k'=1}^{N}\|a^k\|_{L^{2}(\Omega)}.
\|\nabla a^{k'}\|_{L^{\infty}(\Omega)}.
\end{eqnarray*}
From the growth hypothesis on $f$, we deduce that
\begin{eqnarray*}
\Big|\int_{\Omega}f(\sum_{k=1}^{N}d_k^N D(a^k))
:D(a^i)dx \Big| &\leq &  \tilde{c} \|D(a^i)\|_{L^{1}(\Omega)}+c\int_{\Omega}\Big|\sum_{k=1}^{N}d_k^N D(a^k)|^{p-1}|D(a^i)\Big|dx\\ & \leq & C(R,N)\Big(\sum_{k=1}^{N}\|D(a^k)\|_{L^{\infty}(\Omega)}\Big)^{p-1} \|D(a^i)\|_{L^{\infty}(\Omega)}\\ &+&\tilde{c} \|D(a^i)\|_{L^{1}(\Omega)}.\end{eqnarray*}
As $\tilde{\mu}$ is bounded, we estimate in the same way the remaining terms in $\mathcal{F}_{i,j}$ to get
$$|\mathcal{F}_{i,j}(t,c_1^N,...,c_N ^N,d_1^N,...,d_N^N)| \leq C(R,N),$$
where $C(R,N)$ is a constant that does not depend on $t$.\\The standard Caratheodory theory provides the existence of
continuous functions
$(c_1^N,...,c_N^N,d_1^N,...,d_N^N)$ solutions to $(S_N)$ at least for a short time interval with $\frac{d}{dt}(c_1^N,...,c_N^N,d_1^N,...,d_N^N)$ is defined almost everywhere. The uniform estimates that we will derive in the next subsection enable us to extend the solution onto the whole time interval $[0,T]$.
\subsection{Uniform estimates}\label{Uniform}
Multiplying the first equation in $(S_N)$ by $d_i^N(t)$, then taking the sum over $i=1,...,N$, we obtain
$$\frac{1}{2}\frac{d}{dt}\|v^N(t)\|_{L^2(\Omega)}^2+ \int_{\Omega}f(D(v^N)):D(v^N)dx=-\int_{\Omega} \tau^N:D(v^N)dx.$$
The monotonicity hypothesis on $f$ leads to
\begin{equation}
\frac{1}{2}\frac{d}{dt}\|v^N(t)\|_{L^2(\Omega)}^2+\nu (\|D(v^N(t))\|_{L^p(\Omega)}^p+\|D(v^N(t))\|_{L^2(\Omega)}^
 2)\leq \end{equation}
 $$-\int_{\Omega} \tau^N:D(v^N)dx.$$
Multiplying the second equation in $(S_N)$ by $c^N_{j}(t)$, then taking the sum over $j=1,...,N$, we obtain
\begin{eqnarray*}
\frac{1}{2b}\frac{d}{dt}\|\tau^N(t)\|_{L^2(\Omega)}^2+\frac{a}{b} \|\tau^N(t)\|_{L^2(\Omega)}^2=\frac{1}{1-\theta}\int_{\Omega}(\mu(|D(v^N)|^2)-\theta)\tau^N:D(v^N).\end{eqnarray*}
Notice that $(\tau^N. w(v^N)-w(v^N).\tau^N): \tau^N=0$ since $\tau^N$ is symmetric.
Hence, we get
\begin{eqnarray*}
\frac{1}{2}\frac{d}{dt}\|v^N(t)\|_{L^2(\Omega)}^2 &+& \nu  \|D(v^N(t))\|_{L^p(\Omega)}^p+\frac{1}{2b}\frac{d}{dt}\|\tau^N(t)\|_{L^2(\Omega)}^2+\frac{a}{b} \|\tau^N(t)\|_{L^2(\Omega)}^2\\&\leq &\frac{1}{1-\theta}\int_{\Omega}(\mu(|D(v^N)|^2)-1)\tau^N:D(v^N) dx\\ &\leq&  \frac{\lambda}{1-\theta} \|\tau^N(t)\|_{L^2(\Omega)}\|D(v^N(t))\|_{L^2(\Omega)}\\ & \leq &  \frac{\lambda}{1-\theta}\sqrt{\frac{b}{a\nu}}\Big(\sqrt{\nu}\|D(v^N(t))\|_{L^2(\Omega)}\Big).
\Big(\sqrt{\frac{a}{b}}\|\tau^N(t)\|_{L^2(\Omega)}\Big). \end{eqnarray*}
Young inequality implies that
\begin{eqnarray}
\frac{1}{2}\frac{d}{dt}\|v^N(t)\|_{L^2(\Omega)}^2 &+& \nu  \|D(v^N(t))\|_{L^p(\Omega)}^p+\frac{1}{2b}\frac{d}{dt}\|\tau^N(t)\|_{L^2(\Omega)}^2+\frac{a}{b} \|\tau^N(t)\|_{L^2(\Omega)}^2\\ &\leq& \frac{\lambda}{2(1-\theta)}\sqrt{\frac{b}{a\nu}}(\nu \|D(v^N(t))\|_{L^2(\Omega)}^2+ \frac{a}{b}\|\tau^N(t)\|_{L^2(\Omega)}^2).\end{eqnarray}
Let $\gamma:= \frac{\lambda}{2(1-\theta)}\sqrt{\frac{b}{a\nu}}$. Notice that in both cases $i)$ and $ii)$ in Theorem \ref{global1}, we have $\gamma < 1$.
Notice also that to absorb the terms at the right-hand side of $(4.8)$, we do not need that $f$ satisfies the strong monotonicity condition. However, we need just that $f$ satisfies the coercivity condition $f(A):A \geq \nu |A|^p$. In fact, as $\Omega$ is a bounded domain and $p> 2$, then thanks to Holder inequality and Young one
\begin{eqnarray*}\|D(v^N(t))\|_{L^2(\Omega)}^2 &=& \int_{\Omega}|D(v^N(t))|^{2}dx\\ &\leq &
C(|\Omega|)\Big(\int_{\Omega}( |D(v^N(t))|^{2})^{\frac{p}{2}}dx\Big)^{\frac{2}{p}}\\ &\leq & C(|\Omega|)+ \frac{2}{p}\|D(v^N(t))\|_{L^p(\Omega)}^p\\
&\leq & C(|\Omega|)+ \|D(v^N(t))\|_{L^p(\Omega)}^p. \end{eqnarray*}
Finally, we obtain
\begin{eqnarray}\label{energy estimate}
\frac{1}{2}\frac{d}{dt}\|v^N(t)\|_{L^2(\Omega)}^2 &+ &\nu (1-\gamma) \|D(v^N(t))\|_{L^p(\Omega)}^p +\frac{1}{2b}\frac{d}{dt}\|\tau^N(t)\|_{L^2(\Omega)}^2\\
&+& \frac{a}{b}(1-\gamma)\|\tau^N(t)\|_{L^2(\Omega)}^2 \leq C.\end{eqnarray}Particularly, this leads to the fact that$$|(c_1^N,...,c_N^N,d_1^N,...,d_N^N)|^{2}\leq C,$$
where $C$ is a constant that does not depend on $t$ neither on $N$.\\
This uniform boundedness with the continuity of $(c_1^N,...,c_N^N,d_1^N,...,d_N^N)$, we deduce that the functions $(c_1^N,...,c_N^N,d_1^N,...,d_N^N)$ are defined on the whole interval $[0,T]$, see Zeidler \cite{Zeidler} for more details.\\By Korn inequality, we have $$\|\nabla v^N(t)\|_{L^p(\Omega)} \leq C\|D(v^N(t))\|_{L^p(\Omega)},$$
and we deduce that we can extract subsequences such that$$v^N\rightharpoonup^{\ast} v\quad \mbox{weakly in}\,L^{\infty}(0,T,\mathcal{H})$$
$$v^N\rightharpoonup v\quad \mbox{weakly in}\,L^{p}(0,T, \textbf{V}_p)$$
$$\tau^N\rightharpoonup^{\ast} \tau\quad \mbox{weakly in}\quad L^{\infty}(0,T,\widetilde{\mathcal{H}})$$
$$\tau^N\rightharpoonup \tau\quad \mbox{weakly in}\,L^{2}(0,T,\widetilde{\mathcal{H}}).$$
The growth condition on $f$ and the uniform boundedness of the sequence $D(v^N)$ in $L^{p}([0,T]\times \Omega)$ imply that $f(D(v^N))$ is bounded in $L^{p'}([0,T]\times \Omega)$, and thus
$$f(D(v^N)) \rightharpoonup \overline{f(D(v^N))}\quad \mbox{weakly in}\,L^{p'}([0,T]\times \Omega).$$
In addition, as $g(D(v^N))$ is bounded in $L^{p}([0,T]\times \Omega)$, then
$$g(D(v^N)) \rightharpoonup \overline{g(D(v^N))}\quad \mbox{weakly in}\,L^{p}([0,T]\times \Omega).$$
Now, we want to derive some bounds on the pressure. Let us start by the two dimensional case.\\By the Gagliardo Nirenberg inequality, we have
$$\|v^N\|_{L^4(\Omega)}\leq C \|v^N\|_{L^2(\Omega)}^{\frac{1}{2}}\|\nabla v^N\|_{L^2(\Omega)}^{\frac{1}{2}}.$$
Thus, By Holder inequality, we get $$\|v^N\|_{L^4([0,T]\times \Omega)}^4 \leq C \|v^N\|_{L^{\infty}(0,T,L^2(\Omega))}^2 \|\nabla v^N\|_{L^2([0,T]\times \Omega)}^{2}.$$As $p'<2$, then $v^N \otimes v^N$ is bounded in $L^{p'}([0,T]\times \Omega)$.\\
Since $v^N$ is divergence free, we get $$\Delta p^N=  div div\Big(f(D(v^N))-v^N \otimes v^N+ \tau^N\Big).$$
The uniform boundedness of $\tau^{N}$ and $f(D(v^N))$ in $L^{p'}([0,T]\times \Omega)$ imply that $p^{N}$ is also bounded in $L^{p'}([0,T]\times \Omega)$.
Hence, we infer that
$$p^N \rightharpoonup p \quad \mbox{weakly in}\,L^{p'}([0,T]\times \Omega).$$
In the three dimensional case, the Gagliardo-Nirenberg inequality implies
$$\|v^N\|_{L^4(\Omega)}\leq C \|v^N\|_{L^2(\Omega)}^{\frac{1}{4}}\|\nabla v^N\|_{L^2(\Omega)}^{\frac{3}{4}}.$$
Since $p\geq \frac{5}{2}$, then Holder inequality leads to the fact that $v^N \otimes v^N$ is bounded in $L^{p'}([0,T]\times \Omega)$.\\
As above, we get $$p^N \rightharpoonup p \quad \mbox{weakly in}\,L^{p'}([0,T]\times \Omega).$$
Let us remark that at this stage,  in the three dimensional case, we do not really need to have $p \geq \frac{5}{2}$. We need only to have $p> 2$ which lead to the fact that $v^N \otimes v^N$ is bounded in $L^{\frac{4}{3}}(0,T, L^{p'}(\Omega))$ and thus $p^N$ is bounded in $L^{\min(p',\frac{4}{3})}(0,T, L^{p'}(\Omega))$.
\subsection{Strong convergence of $v^N$ in $L^2$}\label{Strong1}
We would like to prove that $(v,\tau,p)$ is still a solution to $(S)$. The difficulties appear when passing to the limit in the nonlinear terms
$$v^N.\nabla v^N;\, v^N.\nabla \tau^N;\,f(D(v^N));\, g(D(v^N)),\, \tau^N. w(v^N)-w(v^N).\tau^N.$$
For the first two terms, we only need to prove that $v^N$ converges strongly to $v$ in $L^{2}([0,T]\times \Omega)$.
The bounds on the pressure and on the term $v^N \otimes v^N$ imply that $\partial_{t} v^N$ is bounded in $L^{r}(0,T, W^{-1, p'}(\Omega))$, for some $r\in ]1,+\infty[$. In addition, $v^N$ is bounded in $L^{2}(0,T, \textbf{V}_2(\Omega))$. As $\textbf{V}_2(\Omega)$ is compactly embedded in $L^2(\Omega)$ which is continuously embedded in $W^{-1, p'}(\Omega)$, thanks to the Aubin-Lions lemma, up to a subsequence $$v^N \rightarrow v\quad \mbox{strongly in}\,L^{2}([0,T]\times \Omega).$$
Hence, we get $$v^N \otimes v^N  \rightarrow v\otimes v \quad \mbox{strongly in}\,L^{1}([0,T]\times \Omega).$$
Since $v^N$ is divergence free, then $$v^N.\nabla \tau^N=\sum_{i=1}^{n}v^N_{i}\partial_{i}\tau^{N}=\sum_{i=1}^{n}\partial_{i}(v^N_{i}\tau^{N})=div(v \tau).$$
As $v^N$ converges strongly to $v$ in $L^{2}([0,T]\times \Omega)$ and $\tau^{N}$ converges weakly to $\tau$ in $L^{2}([0,T]\times \Omega)$, we deduce that $$v^N_{i}\tau^{N} \rightharpoonup v_i \tau \quad \mbox{weakly in}\,L^{1}([0,T]\times \Omega).$$
Thus, it remains to prove that $$\overline{f(D(v^N))}=f(D(v)),\,\,\overline{g(D(v^N))}=g(D(v))$$ $$\mbox {and}$$ $$\overline{\tau^N. w(v^N)-w(v^N).\tau^N}=\tau .w(v)-w(v).\tau.$$
If we prove that $\tau^N$ converges strongly to $\tau$ in $L^{2}([0,T]\times \Omega)$, then our objective will be achieved. To get the strong convergence of $\tau^N$, we will prove first that $\tau^N$ is equi-integrable in $L^2$ uniformly in $t\in [0,T]$:
$$\lim_{M \rightarrow + \infty} \sup_{n}\sup_{t\in [0,T]}\int_{\Omega} |\tau^N|^2 \chi_{\{|\tau^N|\geq M\}}dx=0. $$
\subsection{Equi-integrability of $\tau^N$ in $L^2$}\label{Equi}
We will show the equi-integrability of $\tau^N$ in $L^2$ uniformly in $t\in [0,T]$. The idea is inspired from \cite{Lions+Masm}.\\
Let $R> 0$ be fixed. We decompose $\tau^N$ into the sum of $\psi^{N}$ and $H^N$ that solve respectively the following systems
\begin{eqnarray*}(E^N_1)
\begin{cases}
\partial_t \psi^{N}+v^{N}.\nabla \psi^{N}+a \psi^{N}+ \psi^N. w(v^N)-w(v^N).\psi^N=g(D(v^N)),\\
\psi^{N}(0,x)=\tau_0^N(x)\chi_{\{|\tau_0^N|< R\}},\end{cases}
\end{eqnarray*}
and
\begin{eqnarray*}(E_2^N)
\begin{cases}
\partial_t H^N+v^N.\nabla H^N+a H^N +H^N. w(v^N)-w(v^N).H^N=0\\ H^N(0,x)=\tau_0^N(x)\chi_{\{|\tau_0^N|\geq R\}}.\end{cases}
\end{eqnarray*}
Let $M > 0$, we have
\begin{eqnarray*}
\int_{\Omega}|\tau^N|^{2}\chi_{\{|\tau^N|\geq M\}}dx&=& \int_{\Omega}|\psi^N+H^N|^{2}\chi_{\{|\tau^N|\geq M\}}dx\\
&\leq & 2 \int_{\Omega}|\psi^N|^2 \chi_{\{|\tau^N|\geq M\}}dx+2 \int_{\Omega}|H^N|^2 \chi_{\{|\tau^N|\geq M\}}dx.\end{eqnarray*}
Energy estimates on $H^N$ in the energy space $L^2(\Omega)$ imply that
$$\int_{\Omega} |H^N(t,x)|^2 dx \leq \int_{\Omega} |H^N(0,x)|^2 dx=\|\tau_0^N\chi_{\{|\tau_0^N|\geq R\}}\|_{L^2(\Omega)}^2.$$
Noticing that $\psi^{N}(0,.)\in L^p (\Omega)$, we should have a bound on $\psi^N$ in $L^\infty(0,T, L^p(\Omega))$.\\
We have
\begin{eqnarray*}
\frac{1}{2}\partial_{t}|\psi^{N}|^{2}+a |\psi^{N}|^{2}+\frac{1}{2}v^N .\nabla |\psi^{N}|^{2}=g(D(v^N)):\psi^N.\end{eqnarray*}
Multiplying this equation by $p |\psi^N|^{p-2}$, we get
$$\partial_t |\psi^N|^{p}+a p |\psi^N|^{p}+v^N. \nabla |\psi^N|^{p}=p |\psi^N|^{p-2}g(D(v^N)):\psi^N.$$
As $\tilde{\mu}$ is bounded, integrating over $\Omega$, we deduce that
$$\partial_t\int_{\Omega}|\psi^N(t,x)|^{p}dx\leq p \int_{\Omega} |\psi^N(t,x)|^{p-1}|D(v^N)(t,x)|dx.$$
Thanks to Holder inequality, we get
$$\partial_t\int_{\Omega}|\psi^N(t,x)|^{p}dx\leq p \Big(\int_{\Omega} |\psi^N (t,x)|^{p}dx\Big)^{\frac{p-1}{p}}\Big(\int_{\Omega}|D(v^N)(t,x)|^p dx\Big)^{\frac{1}{p}}.$$
Multiplying this equation by $\Big(\int_{\Omega}|\psi^N(t,x)|^{p}dx\Big)^{\frac{1}{p}-1}$, we obtain$$\partial_t \|\psi^N(t,.)\|_{L^p(\Omega)}\leq
\Big(\int_{\Omega}|D(v^N)(t,x)|^p dx\Big)^{\frac{1}{p}},$$
and thus by Holder inequality in time $$\|\psi^N(t,.)\|_{L^p(\Omega)}\leq \|\psi^N(0,.)\|_{L^p(\Omega)}+ T^{1-\frac{1}{p}}\|D(v^N)\|_{L^p([0,T]\times \Omega)} .$$
As $D(v^N)$ is uniformly bounded in $L^p([0,T]\times \Omega)$, we get $$\|\psi^N(t,.)\|_{L^p(\Omega)}\leq \|\psi^N(0,.)\|_{L^p(\Omega)}+C.$$
Applying Holder inequality, we obtain
\begin{eqnarray*}\int_{\Omega}|\psi^N|^2 \chi_{\{|\tau^N|\geq M\}}dx &\leq& \Big(\int_{\Omega}|\psi^N(t,x)|^{p}dx \Big)^{\frac{2}{p}} \Big(\int_{\Omega} \chi_{\{|\tau^N|\geq M\}}dx\Big)^{1-\frac{2}{p}}\\
&\leq&   \|\psi^N\|_{L^\infty(0,T, L^p(\Omega))}^2\Big(\frac{1}{M^2}  \int_{\Omega}|\tau^{N}(t,x)|^2 dx \Big)^{1-\frac{2}{p}}\\
&\leq& C(\|\psi^N(0,.)\|_{L^p(\Omega)}^2+C)\frac{1}{M^{2(1-\frac{2}{p})}}\|\tau^N\|_{L^\infty(0,T,L^2(\Omega))}^{2(1-\frac{2}{p})}.\end{eqnarray*}
Since $\tau^N$ is uniformly bounded in $L^\infty(0,T,L^2(\Omega))$, then we obtain
\begin{eqnarray*}\int_{\Omega}|\psi^N|^2 \chi_{\{|\tau^N|\geq M\}}dx \leq C \frac{1}{M^{2(1-\frac{2}{p})}}(\|\psi^N(0,.)\|_{L^p(\Omega)}^2+C).\end{eqnarray*}
But, we have as $p> 2$\begin{eqnarray*}
\|\psi^N(0,.)\|_{L^p(\Omega)}^p&=&\int_{\Omega} |\tau^{N}(0,x)|^p \chi_{\{|\tau_0^N|< R\}}dx\\
&\leq& \int_{\Omega} |\tau^{N}(0,x)|^{p-2}|\tau^{N}(0,x)|^{2}\chi_{\{|\tau_0^N|< R\}}dx\\
&\leq& R^{p-2}\|\tau^{N}(0,.)\|_{L^2(\Omega)}^2.\end{eqnarray*}
Since, we have $\| P^N\|_{\mathcal{L}(\widetilde{\mathcal{H}},\widetilde{\mathcal{H}})} \leq 1$, we deduce that
$$\|\psi^N(0,.)\|_{L^p(\Omega)}\leq R^{\frac{p-2}{p}}\|\tau^{N}(0,.)\|_{L^2(\Omega)}^{\frac{2}{p}}\leq R^{\frac{p-2}{p}}\|\tau_0\|_{L^2(\Omega)}^{\frac{2}{p}}.$$
Finally, for a fixed $R> 0$, we get
$$\sup_N \sup_{t\in [0,T]} \int_{\Omega}|\tau^N|^{2}\chi_{\{|\tau^N|\geq M\}}dx \leq 2 \sup_N\|\tau_0^N\chi_{\{|\tau_0^N|\geq R\}}\|_{L^2(\Omega)}^2+ C \frac{1}{M^{2(1-\frac{2}{p})}}( R^{\frac{2(p-2)}{p}} \|\tau_0\|_{L^2(\Omega)}^{\frac{4}{p}}  +C) .$$
As $p> 2$, we obtain
\begin{equation}\label{ine1}
\limsup_{M} \sup_N\sup_{t\in [0,T]}\int_{\Omega}|\tau^N|^{2}\chi_{\{|\tau^N|\geq M\}}dx \leq 2 \sup_N\|\tau_0^N\chi_{\{|\tau_0^N|\geq R\}}\|_{L^2(\Omega)}^2,\end{equation}
which hold true for every $R \in ]0,+\infty[$.
The fact that $\tau_0^N$ converges strongly to $\tau_0$ in $L^2(\Omega)$ implies that $|\tau_0^N|^2$ converges strongly in $L^1(\Omega)$ and thus weakly in $L^1(\Omega)$, and thus $|\tau_0^N|^2$ is equi-integrable in $L^1(\Omega)$.\\
Consequently, let $R\rightarrow +\infty$ in inequality \eqref{ine1}, we infer that
$$0\leq \liminf_{M} \sup_N\sup_{t\in [0,T]}\int_{\Omega}|\tau^N|^{2}\chi_{\{|\tau^N|\geq M\}}dx\leq \limsup_{M} \sup_N\sup_{t\in [0,T]}\int_{\Omega}|\tau^N|^{2}\chi_{\{|\tau^N|\geq M\}}dx\leq 0.$$
This means that $\tau^N$ is equi-integrable in $L^2(\Omega)$ uniformly in $t\in [0,T]$.
\subsection{Strong convergence of $\tau^N$ in $L^2$}\label{Strong2}
We will first focus on the system $(S)$.\\As $\tau^N$ is equi-integrable in $L^2(\Omega)$ uniformly in $t\in [0,T]$, then $|\tau^N-\tau|^2$ converges weakly in $L^1([0,T]\times \Omega)$. So let $\eta:=\overline{|\tau^N-\tau|^2}\in L^\infty(0,T, L^1(\Omega))$.\\Notice that if $\eta=0$, then $$|\tau^N-\tau|^2 \rightarrow 0 \quad\mbox{strongly in}\quad L^1([0,T]\times \Omega).$$Hence, our aim will be to show that $\eta=0$.\\
Multiplying by $\tau^N$ the equation satisfied by $\tau^N$, we get
\begin{equation}\label{t1}(\partial_{t}+v^N.\nabla) |\tau^{N}|^{2}+ 2a |\tau^{N}|^{2}=2g(D(v^N)):\tau^N.\end{equation}
Notice here that the term $v^N. \nabla \tau^N :\tau^N$ has a sense since $\nabla \tau^N \in L^\infty, v^N \in L^2, \tau^N \in L^2$  for a fixe $N \in \mathbf{N}$.\\
Let us introduce the unique a.e flow in the sense of DiPerna and Lions \cite{dip+lio} of $v^N$, solution of $$\partial_{t} X^N(t,x)= v^N(t,X(t,x)), \quad X^N(0,x)=x.$$
We also denote by $X$ the a.e flow of $v$.\\
Equation \eqref{t1} can be written as \begin{eqnarray*}\partial_t [|\tau^{N}|^{2}(t,X^N(t,x))]+ 2a [|\tau^{N}|^{2}(t,X^N(t,x))]= 2\Big(g(D(v^N)):\tau^N\Big)
(t,X^N(t,x)).\end{eqnarray*}
Passing to the limit weakly, we get
\begin{eqnarray*}\partial_t \overline{|\tau^{N}|^{2}(t,X^N(t,x))}+ 2a \overline{|\tau^{N}|^{2}(t,X^N(t,x))}=2\overline{\Big(g(D(v^N)):\tau^N\Big)
(t,X^N(t,x))}.\end{eqnarray*}
The stability of the notion of a.e flow with respect to the weak limit of $v^N$ to $v$ implies that $X^N(t,x)$ converges to $X(t,x)$ in $L^1_{loc}$ and also that $(X^N(t)^{-1})(x)$ converges to $(X(t)^{-1})(x)$ in $L^1_{loc}$. Thus, we get
$$\overline{|\tau^{N}|^{2}(t,X^N(t,x))}=\overline{|\tau^{N}|^{2}(t,X(t,x))}=(|\tau|^2+\eta)(t,X(t,x)).$$
Finally, we obtain
\begin{equation}\label{t2} \frac{1}{2}(\partial_{t}+v.\nabla)(|\tau|^{2}+\eta)+a(|\tau|^{2}+\eta)=
\overline{g(D(v^N)):\tau^N}.\end{equation}
In addition, since $v^N.\nabla \tau^N$ converges in $D'([0,T]\times \Omega)$ to $v.\nabla \tau$, the equation satisfied by $\tau$ is then $$\partial_t \tau +a \tau+ v.\nabla \tau +\overline{\tau^N.w(v^N)-w(v^N).\tau^N}=\overline{g(D(v^N))}, $$
which implies that \begin{equation}\label{t3} \frac{1}{2}(\partial_{t}+v.\nabla)|\tau|^{2}+a |\tau|^{2}+\overline{\tau^N.w(v^N)-w(v^N).\tau^N}:\tau=\overline{g(D(v^N))}:\tau.\end{equation}
Subtracting \eqref{t3} from \eqref{t2} gives the equation satisfied by $\eta$
\begin{equation}\label{t4}
\partial_{t}\eta +2 a \eta +div (\eta v) -2 \overline{\tau^N.w(v^N)-w(v^N).\tau^N}:\tau=2\Big(\overline{g(D(v^N)):\tau^N}-\overline{g(D(v^N))}:\tau\Big).
\end{equation}
Notice that the term $\eta v$ makes sense in the sense of distributions in the two dimensional case as $W^{1,p}(\Omega)\subset L^\infty(\Omega)$ and in the three dimensional case for $p> 3$, and also the term $\overline{\tau^N.w(v^N)-w(v^N).\tau^N}:\tau$ is not defined in the sense of distributions. To overcome this difficulty for any $p$, we use a renormalized form of \eqref{t4} by multiplying such equation by $\frac{1}{(1+\eta)^2}$, namely
\begin{eqnarray*}(\partial_{t}+v.\nabla)\zeta +\frac{2a}{1+\eta} \zeta  &=&\frac{2}{(1+\eta)^2}\Big(\overline{g(D(v^N)):\tau^N}-\overline{g(D(v^N))}:\tau\\ &+& \overline{\tau^N.w(v^N)-w(v^N).\tau^N}:\tau \Big),\end{eqnarray*} where $\zeta=\frac{\eta}{1+\eta}\in L^\infty ([0,T]\times \Omega)$.\\
Let us remark that$$T_1:=\overline{g(D(v^N)):\tau^N}-\overline{g(D(v^N))}:\tau=\overline{\Big(g(D(v^N))-
g(D(v))\Big):\Big(\tau^N-\tau\Big)},$$ and $$T_2:=\overline{\tau^N.w(v^N)-w(v^N).\tau^N}:\tau=\overline{(\tau^N-\tau).(w(v^N)-w(v))-(w(v^N)-w(v)).(\tau^N-\tau)}:\tau.$$
As $\tau^N$ is equi-integrable in $L^2([0,T]\times \Omega)$ and $w(v^N)$ is bounded in $L^2([0,T]\times \Omega)$, then $(\tau^N-\tau).(w(v^N)-w(v))$ is
equi-integrable in $L^1([0,T]\times \Omega)$ and thus converges weakly in $L^1([0,T]\times \Omega)$. Let $\nu_{t,x}$ be the Young measure generated by $(\tau^N,w(v^N))$, then we have
\begin{eqnarray*}
\overline{(\tau^N-\tau).(w(v^N)-w(v))}=\int_{M^{n\times n}\times M^{n\times n}}(\alpha-\tau).(\beta-w(v))d \nu_{t,x}(\alpha,\beta)
\end{eqnarray*}
By Cauchy-Schwartz inequality, we get
\begin{eqnarray*}
|\overline{(\tau^N-\tau).(w(v^N)-w(v))}|&\leq &\Big(\int_{M^{n\times n}\times M^{n\times n}}|\alpha-\tau|^2 d \nu_{t,x}(\alpha,\beta)\Big)^{\frac{1}{2}}\\
&\times& \Big(\int_{M^{n\times n}\times M^{n\times n}}|\beta-w(v)|^2 d \nu_{t,x}(\alpha,\beta)\Big)^{\frac{1}{2}}.
\end{eqnarray*}
The equi-integrability of $\tau^N$ in $L^2([0,T]\times \Omega)$ implies that
$$\eta = \int_{M^{n\times n}\times M^{n\times n}}|\alpha-\tau|^2 d \nu_{t,x}(\alpha,\beta).$$
In addition, $|w(v^N)|^2$ is equi-integrable in $L^1([0,T]\times \Omega)$ since it is bounded in $L^{\frac{p}{2}}([0,T]\times \Omega)$ with $p> 2$, thus
$$\overline{|w(v^N)-w(v)|^2}=\int_{M^{n\times n}\times M^{n\times n}}|\beta-w(v)|^2 d \nu_{t,x}(\alpha,\beta).$$
Consequently, we obtain
\begin{eqnarray*}
|\overline{(\tau^N-\tau).(w(v^N)-w(v))}|\leq \sqrt{\eta} \overline{|w(v^N)-w(v)|^2}^{\frac{1}{2}},
\end{eqnarray*}
and we remark that $$\frac{1}{(1+\eta)^2}|\overline{(\tau^N-\tau).(w(v^N)-w(v))}:\tau|\leq \frac{\sqrt{\eta}}{(1+\eta)^2} |\tau|\overline{|w(v^N)-w(v)|^2}^{\frac{1}{2}},$$
with $\frac{\sqrt{\eta}}{(1+\eta)^2} \in L^\infty([0,T]\times \Omega)$, $|\tau|\in L^2([0,T]\times \Omega)$ and $\overline{|w(v^N)-w(v)|^2}^{\frac{1}{2}}\in L^2([0,T]\times \Omega)$.\\
As $\tau^N$ is equi-integrable in $L^2([0,T]\times \Omega)$ and $D(v^N)$ is bounded in $L^2([0,T]\times \Omega)$, then $(\tau^N-\tau).(g(D(v^N))-g(D(v)))$ is
equi-integrable in $L^1([0,T]\times \Omega)$. Let $\tilde{\nu}_{t,x}$ be the Young measure generated by $(\tau^N,D(v^N))$, then we have
$$T_1=\int_{M^{n\times n}\times M^{n\times n}}(\alpha-\tau).(g(\beta)-g(D(v)))d \tilde{\nu}_{t,x}(\alpha,\beta).$$
As the derivative of $g$ is bounded, then $g$ is Lipschitz and the same calculus as for the term $T_2$ imply that
$$T_1 \leq C \sqrt{\eta} \overline{|D(v^N)-D(v)|^2}^{\frac{1}{2}}.$$
Finally, we get
\begin{eqnarray*}
(\partial_{t}+v.\nabla)\zeta \leq C \frac{\sqrt{\eta}}{(1+\eta)^2} \Big(\overline{|D(v^N)-D(v)|^2}^{\frac{1}{2}}+  |\tau|\overline{|w(v^N)-w(v)|^2}^{\frac{1}{2}}\Big).
\end{eqnarray*}
Now, we should estimate $\overline{|D(v^N)-D(v)|^2}$ and $\overline{|w(v^N)-w(v)|^2}$ in terms of $\eta$ with the help of the first equation in $(S)$. More precisely, we will prove that $$\overline{(f(D(v^{N}))-f(D(v))):(D(v)-D(v^{N}))}=\overline{(\tau^{N}-\tau):(D(v^{N})-D(v))}.$$
Multiplying the first equation in $(S)$ by $v^N$ and using the identity $$div(A)u=div(Au)-A:D(u),\, A \in M^{n\times n}, u \in \mathbf{R}^n,$$ we get
\begin{equation} \label{first}\frac{1}{2}\partial_t |v^N|^{2}-div(f(D(v^{N})) v^N)+ f(D(v^{N})):D(v^N)+ v^N .\nabla v^N .v^N+div(p^N v^N)=div( \tau^N v^N)-\tau^N:D(v^N).
\end{equation}
The strong convergence of $v^N$ to $v$ in $L^2([0,T]\times \Omega)$ implies that $|v^N|^2$ converges strongly in $L^1([0,T]\times \Omega)$ to $|v|^2$.\\In addition, thanks to the weak convergence of $\tau^N$ to $\tau$ in $L^2([0,T]\times \Omega)$, we can conclude that $$ \tau^N v^N \rightharpoonup \tau v \quad \mbox{weakly in}\,L^{1}([0,T]\times \Omega).$$
But, $v^N$ is uniformly bounded in $L^p(0,T, \textbf{V}_p)$ and $\partial_t v^N$  is uniformly bounded in $L^{p'}(0,T, W^{-1,p'}(\Omega))$ in the two dimensional case and $L^{\min{(p',\frac{4}{3})}}(0,T, W^{-1,p'}(\Omega))$ in the three dimensional case. Since $\textbf{V}_p$ is compactly embedded into $L^p$ which is continuously embedded into $W^{-1,p'}(\Omega)$, Rellich-Kondrachov Theorem implies that $$v^N \rightarrow  v\quad \mbox{strongly in}\,L^{p}([0,T]\times \Omega).$$
As $f(D(v^n))$ is uniformly bounded in $L^{p'}([0,T]\times \Omega)$, we obtain
$$ f(D(v^n)) v^N \rightharpoonup \overline{f(D(v^n))}v \quad \mbox{weakly in}\,L^{1}([0,T]\times \Omega).$$
In the two dimensional case, $p^N$ and $v^N \otimes v^N$ are uniformly bounded in $L^{p'}([0,T]\times \Omega)$. Thus, one gets
$$ p^N v^N \rightharpoonup p v \quad \mbox{weakly in}\,L^{1}([0,T]\times \Omega).$$
Notice that thanks to the divergence free condition on $v^N$, $$v^N .\nabla v^N .v^N=\frac{1}{2}div(|v^N|^2 v^N),$$ and thus
$$ |v^N|^2 v^N \rightharpoonup \overline{|v^N|^2} v=|v|^2 v \quad \mbox{weakly in}\,L^1([0,T]\times \Omega).$$
In the three dimensional case, for $p \geq \frac{5}{2}$, $p^N$ and $v^N \otimes v^N$ are uniformly bounded in $L^{p'}([0,T]\times \Omega)$.
Notice that at this stage appears the bound on $p$ when the dimension is three. In fact, we have
\begin{eqnarray*}
\|v^N \otimes v^N\|^{p'}_{L^{p'}(\Omega)} &\leq& C \|v^N \|^{2p'}_{L^{4}(\Omega)}\\ &\leq&\|v^N \|^{\frac{p'}{2}}_{L^{2}(\Omega)}
\| \nabla v^N \|^{\frac{3p'}{2}}_{L^{2}(\Omega)}\\ &\leq& C \|v^N \|^{\frac{p'}{2}}_{L^{2}(\Omega)} \| \nabla v^N \|^{\frac{3p'}{2}}_{L^{p}(\Omega)}.\end{eqnarray*}
Thus, it suffices to have $\frac{3p'}{2} \leq p$ that is $p \geq \frac{5}{2}$.\\Taking the weak limit of \eqref{first}, we obtain
\begin{equation} \label{second}\frac{1}{2}\partial_t |v|^{2}-div(\overline{f(D(v^{N}))} v)+ \overline{f(D(v^{N})):D(v^N)}+ v .\nabla v .v+div(p v)=div( \tau v)-\overline{\tau^N:D(v^N)}.\end{equation}
On the other hand, we take the weak limit of the first equation in $(S)$ knowing that $v^N.\nabla v^N $ converges to $v.\nabla v$
$$\partial_t v +v.\nabla v-div (\overline{f(D(v^{N}))}) +\nabla p=div \tau,$$
which is multiplied by $v$ leads to
\begin{equation}\label{third}\frac{1}{2}\partial_t |v|^{2}-div(\overline{f(D(v^{N}))} v)+\overline{f(D(v^{N}))}:D(v)+ v.\nabla v.v++div(p v)=div( \tau v)-\tau:D(v).\end{equation}
Subtracting \eqref{third} from \eqref{second}, we get
$$\overline{f(D(v^{N})):D(v^N)}-\overline{f(D(v^{N}))}:D(v)=-\overline{\tau^N:D(v^N)}+\tau:D(v).$$
Or equivalently
$$\overline{(f(D(v^{N}))-f(D(v))):(D(v)-D(v^{N}))}=\overline{(\tau^{N}-\tau):(D(v^{N})-D(v))}.$$
For the system $(S_2)$, under the hypotheses on $f$ in Theorem \ref{global2}, since $f$ is monotone, we deduce that $$\overline{(\tau^{N}-\tau):(D(v^{N})-D(v))}\leq 0,$$
and thus, the equation satisfied by $\zeta$ is reduced to $$(\partial_t +v. \nabla) \zeta \leq\frac{2}{(1+\eta)^2} \overline{(\tau^{N}-\tau):(D(v^{N})-D(v))}\leq 0, \, \zeta(0,.)=0, $$ which implies that $\zeta= \eta=0$.\\
Now, for the system $(S_1)$, under the strong monotonicity condition on $f$ in Theorem \ref{global1}, we deduce that
\begin{eqnarray*}
\overline{|D(v^N)-D(v)|^2} \leq |\overline{(\tau^{N}-\tau):(D(v^{N})-D(v))}| &\leq &\sqrt{\eta}\overline{|D(v^N)-D(v)|^2}^{\frac{1}{2}}\\
&\leq& \frac{1}{2}\eta + \frac{1}{2} \overline{|D(v^N)-D(v)|^2}.
\end{eqnarray*}
Thus, we get
\begin{eqnarray*}
\overline{|D(v^N)-D(v)|^2} \leq \eta,
\end{eqnarray*}
and
$$(\partial_t +v.\nabla )\zeta \leq C \zeta,$$
which implies thanks to Gronwall's lemma and to the fact that $\zeta(0,.)=0$ that $\zeta =\eta=0$.\\
Notice that if the quadratic term $\tau^N. w(v^N)-w(v^N).\tau^N$ is present, then the difficulty is how to estimate $\overline{|w(v^N)-w(v)|^2}$ in terms of $\eta$.
\subsection{Weak limits of non linear terms involving $D(v^N)$}\label{Weak}
Under the hypotheses on $f$ in Theorem \ref{global1}, we have $$\overline{|D(v^N)-D(v)|^2}\leq \eta =0.$$
As $|D(v^N)-D(v)|^2$ is bounded in $L^{\frac{p}{2}}([0,T]\times \Omega)$ with $p> 2$, then it is equi-integrable in $L^{1}([0,T]\times \Omega)$, and hence
$$D(v^N) \rightarrow D(v)\quad \mbox{strongly  in}\,L^{2}([0,T]\times \Omega).$$
Vitali's lemma implies that $$\overline{f(D(v^N))}=f(D(v)),\, \overline{g(D(v^N))}=g(D(v)).$$
Let us now focus on the system $(S_2)$. Under the hypotheses on $f$ in Theorem \ref{global2}, we have two ways to get $\overline{f(D(v^N))}$.\\
Let us begin by the simplest way. Let $\mu_{t,x}$ be the Young measure generated by the sequence $D(v^N)$ and let $G(\lambda, t,x):=(f(\lambda)-f(D(v)(t,x))):(\lambda -D(v)(t,x))\geq 0$.\\
Remark that $$\int_{M^{n\times n}} G(\lambda, t,x)\mu_{t,x}(\lambda)=\int_{M^{n\times n}}\liminf_{\delta  \rightarrow 0} \frac{G(\lambda, t,x)}{1+ \delta G(\lambda, t,x)}d\mu_{t,x}(\lambda).$$
 By Fatou's lemma, we get $$\int_{M^{n\times n}} G(\lambda, t,x)d\mu_{t,x}(\lambda) \leq \liminf_{\delta  \rightarrow 0} \int_{M^{n\times n}}\frac{G(\lambda, t,x)}{1+ \delta G(\lambda, t,x)}d\mu_{t,x}(\lambda).$$
Notice that$$\frac{G(D(v^N), t,x)}{1+ \delta G(D(v^N), t,x)} \leq \frac{1}{\delta}\frac{G(D(v^N), t,x)}{1+ G(D(v^N), t,x)} \leq\frac{1}{\delta}, \, \delta \in ]0,1].$$
Hence, for a fixed $\delta \in ]0,1]$, the sequence $\frac{G(D(v^N), t,x)}{1+ \delta G(D(v^N), t,x)} $ is bounded in $L^\infty ([0,T]\times \Omega)$ and thus equi-integrable in $L^1 ([0,T]\times \Omega)$.\\ We can claim that $$\overline{\frac{G(D(v^N), t,x)}{1+ \delta G(D(v^N), t,x)} }=\int_{M^{n\times n}}
\frac{G(\lambda, t,x)}{1+ \delta G(\lambda, t,x)}d\mu_{t,x}(\lambda), \, \forall \delta \in ]0,1].$$
Consequently, we get
$$\int_{M^{n\times n}} G(\lambda, t,x)d\mu_{t,x}(\lambda) \leq \liminf_{\delta  \rightarrow 0} \overline{\frac{G(D(v^N), t,x)}{1+ \delta G(D(v^N), t,x)} }\leq \overline{G(D(v^N), t,x)}.$$
But, thanks to the strong convergence of $\tau^N$ to $\tau$ in $L^2([0,T]\times \Omega)$ and to the weak convergence of $D(v^N)$ to $D(v)$ in $L^2([0,T]\times \Omega)$, we deduce that $$\int_{M^{n\times n}} G(\lambda, t,x)d\mu_{t,x}(\lambda) \leq \overline{G(D(v^N), t,x)}=-\overline{(\tau^N-\tau)(D(v^N)-D(v))}=0.$$ This means that the div-curl inequality used in \cite{Drey+Hung} and \cite{Drey+Hung'}, which is the key ingredient to handle the limits of non-linear terms, remains true in our case and enables us to get
$\overline{f(D(v^N))}= f(D(v))$.\\
The second way is to proceed as in \cite{Drey+Hung} and \cite{Drey+Hung'} to get the div-curl inequality.\\
As $f$ is monotone, then the negative part of $(f(D(v^N))-f(D(v))):(D(v^N)-D(v))$ is null and thus weakly relatively compact in $L^1$. Applying Theorem \ref{fatou}, we infer that
\begin{eqnarray*}\int_{0} ^{t}\int_{\Omega}\int_{M^{n\times n}} (f(\lambda)-f(D(v)))&:&(\lambda-D(v))d\mu_{s,x}(\lambda)dx ds \\ &\leq& \liminf_{n \rightarrow +\infty}\int_{0} ^{t}\int_{\Omega} (f(D(v^N))-f(D(v))):(D(v^N)-D(v))dx ds\\ &\leq & \limsup_{n \rightarrow +\infty}\int_{0} ^{t}\int_{\Omega} (f(D(v^N))-f(D(v))):(D(v^N)-D(v)) dx ds \\ &\leq &
- \liminf_{n \rightarrow +\infty}\int_{0} ^{t}\int_{\Omega}f(D(v)):D(v^N)dx ds+ \int_{0} ^{t}\int_{\Omega}f(D(v)):D(v)dxds\\
&+&\limsup_{n \rightarrow +\infty}\int_{0} ^{t}\int_{\Omega}f(D(v^N)):D(v^N)dx ds\\&-&
\liminf_{n \rightarrow +\infty}\int_{0} ^{t}\int_{\Omega}f(D(v^N)):D(v)dxds.\\\end{eqnarray*}
The weak convergence in $L^p$ of $D(v^N)$ to $D(v)$ and the fact that $f(D(v))$ belongs to $L^{p'}$ imply that
$$\liminf_{n \rightarrow +\infty}\int_{0} ^{t}\int_{\Omega}f(D(v)):D(v^N)dxds=\int_{0} ^{t}\int_{\Omega}f(D(v)):D(v)dx ds.$$
Moreover, we have $$\liminf_{n \rightarrow +\infty}\int_{0} ^{t}\int_{\Omega}f(D(v^N)):D(v)dxds=\int_{0} ^{t}\int_{\Omega}\overline{f(D(v^N))}:D(v)dxds.$$
Testing the first equation in $(S)$ by $v^N$ and integrating over $[0,t]\times \Omega$, we get
\begin{eqnarray*}\int_{0} ^{t}\int_{\Omega}f(D(v^N)):D(v^N)dx ds &=&\int_{0} ^{t}\int_{\Omega}(-\frac{1}{2}\partial_t |v^N|^2-\tau^N :D(v^N))dx ds\\ &\leq & -\int_{0} ^{t}\int_{\Omega}\tau^N :D(v^N)dx ds+\frac{1}{2}\|v^N_0\|_{L^2(\Omega)}^2- \frac{1}{2}\|v^N(t,.)\|_{L^2(\Omega)}^2.\end{eqnarray*}
The same arguments as for the classical Navier-Stokes equation lead to $$v^N (t,.) \rightharpoonup v(t,.) \,\,\mbox{weakly in} \,\,L^2 , \forall t \in [0,T].$$
Thanks to the lower semi-continuity of the norm, we have $$\|v(t,.)\|_{L^2(\Omega)}\leq \liminf_{n \rightarrow +\infty }\|v^N(t,.)\|_{L^2(\Omega)}.$$
Using the fact that $v^N_0$ converges strongly to $v_0$ in $L^2$, we obtain
\begin{eqnarray*}\limsup_{n \rightarrow +\infty}\int_{0} ^{t}\int_{\Omega}f(D(v^N)):D(v^N)dxds&\leq & \frac{1}{2}\|v_0\|_{L^2(\Omega)}^2 +\limsup_{n \rightarrow +\infty}
\int_{0} ^{t}\int_{\Omega}-\tau^N :D(v^N)dx ds \\ &-& \frac{1}{2}\|v(t,.)\|_{L^2(\Omega)}^2.\end{eqnarray*}
Knowing that $$\frac{1}{2}\partial_t |v|^{2}-div(\overline{f(D(v^{N}))} v)+\overline{f(D(v^{N}))}:D(v)+ v.\nabla v.v+div(p v)=div( \tau v)-\tau:D(v),$$
we get by integration over $[0,t]\times \Omega$
$$\int_{0} ^{t}\int_{\Omega}\overline{f(D(v^{N}))}:D(v)dx ds =-\frac{1}{2}\|v(t,.)\|_{L^2(\Omega)}^2+ \frac{1}{2}\|v_0\|_{L^2(\Omega)}^2-\int_{0} ^{t}\int_{\Omega} \tau:D(v)dx ds.$$
Finally, we get
\begin{eqnarray*}\int_{0} ^{t}\int_{\Omega}\int_{M^{n\times n}} (f(\lambda)-f(D(v)))&:&(\lambda-D(v))d\mu_{s,x}(\lambda)dx ds \\ &\leq& \liminf_{n \rightarrow +\infty}\int_{0} ^{t}\int_{\Omega} (f(D(v^N))-f(D(v))):(D(v^N)-D(v))dx ds\\ &\leq & \limsup_{n \rightarrow +\infty}
\int_{0} ^{t}\int_{\Omega}-\tau^N :D(v^N)dx ds+ \int_{0} ^{t}\int_{\Omega} \tau:D(v)dx ds.\end{eqnarray*}
Since $\tau^N$ converges strongly to $\tau$ in $L^2$ and $D(v^N)$ converges weakly in $L^2$ to $D(v)$, then $\tau^N :D(v^N)$ converges weakly in $L^1$ to $\tau: D(v)$, particularly $$\lim_{n \rightarrow +\infty}\int_{0} ^{t}\int_{\Omega}\tau^N :D(v^N)dx ds=\limsup_{n \rightarrow +\infty}
\int_{0} ^{t}\int_{\Omega}\tau^N :D(v^N)dx ds=\int_{0} ^{t}\int_{\Omega}\tau :D(v)dx ds.$$
Hence
\begin{eqnarray*}\int_{0} ^{t}\int_{\Omega}\int_{M^{n\times n}} (f(\lambda)-f(D(v)))&:&(\lambda-D(v))d\mu_{s,x}(\lambda)dx ds \\ &\leq& \liminf_{n \rightarrow +\infty}\int_{0} ^{t}\int_{\Omega} (f(D(v^N))-f(D(v))):(D(v^N)-D(v))dx ds\\ &\leq & 0.\end{eqnarray*}
This ends the proofs of Theorems \ref{global1}-\ref{global2}.


\vspace{1cm}

\noindent{\bf Acknowledgments.} {\tt  We are very grateful  to Professor Nader Masmoudi for  interesting
 discussions around the questions dealt with in this paper.}


\end{document}